\documentclass[letter, 11pt, psamsfonts]{amsart}
\usepackage{color}

\usepackage{titletoc}

\usepackage{amsmath}
\usepackage{amsthm}
\usepackage{amssymb}
\usepackage{amscd}
\usepackage{amsfonts}
\usepackage{amsbsy}
\usepackage{graphicx}
\usepackage{calc}
\usepackage{setspace}
\usepackage{bbm}
\usepackage{paralist}
\usepackage{setspace}
\usepackage{verbatim}
\usepackage{euscript}
\usepackage{mathrsfs}

\newcommand{\R}{\ensuremath{\mathbb{R}}}
\newcommand{\C}{\ensuremath{\mathbb{C}}}
\newcommand{\N}{\ensuremath{\mathbb{N}}}
\newcommand{\Q}{\ensuremath{\mathbb{Q}}}

\newcommand{\p}{\partial}

\def\B{{\mathscr{B}}}

\newtheorem {theorem} {Theorem}%[section]

\newcommand{\fb}{{\bf f}}

\newcommand{\vv}{{\bf V}}
\newcommand{\xb}{{\bf x}}

\newcommand{\bP}{\mathbb P}
\newcommand{\bC}{\mathbb C}
\newcommand{\Var}{\mathcal V}
\newcommand{\sL}{\mathcal L}
\newcommand{\sK}{\mathcal K}
\newcommand{\sH}{\mathcal H}
\usepackage[colorlinks=true, linkcolor=blue]{hyperref}

\title[Symbolic-numerical approach to the center-focus problem]{A hybrid symbolic-numerical approach\\ 
to the center-focus problem}
\author[A. Mahdi, C. Pessoa, and J.D. Hauenstein]{Adam Mahdi$^1$, Claudio Pessoa$^2$ and Jonathan D. Hauenstein$^3$}
\address{$^1$ Institute of Biomedical Engineering, University of Oxford, UK;  and  Faculty of Applied Mathematics, AGH University of Science of Technology, Poland} \email{adam.mahdi@eng.ox.ac.uk} 

\address{$^2$ Universidade Estadual Paulista,
Departamento de Matem\'atica,  IBILCE/UNESP, 
Rua Cristov\~ao Colombo, 2265, 15.054-–000, S\~ao Jos\'e do Rio Preto, SP, Brazil} \email{pessoa@ibilce.unesp.br} 

\address{$^3$ Department of Applied and Computational Mathematics and Statistics, University of Notre Dame, Notre Dame, IN 46556 USA} \email{hauenstein@nd.edu} 

\subjclass[2010]{34C05
34A34} \keywords{center-focus problem, center manifold, first
integral, numerical algebraic geometry}
\date{\today}

\begin{document}

\maketitle

\onehalfspacing
%\tableofcontents

\begin{abstract}
We propose a new hybrid symbolic-numerical approach to the center-focus problem.  The method allowed us to obtain center conditions for a three-dimensional system of differential equations, which was previously not possible using traditional, purely symbolic computational techniques. 
\end{abstract}

%%%%%%%%%%%%%%%%%%%%%%%%%%%%%%%%%%%
\section{Introduction}\label{sec:01}
%%%%%%%%%%%%%%%%%%%%%%%%%%%%%%%%%%%

%\subsection{Stability for nonhyperbolic singularity}
%\subsection{Computational challanges}
%\subsection{Symbolic computation}
%\subsection{Numerical algebraic geometry}

\subsection{Background}
Determination of the local stability of an isolated singular point for a system of ordinary differential equations (ODEs) is one of the fundamental problems encountered across  various branches of applied sciences and engineering. For a system
\begin{equation}\label{sys}
\dot \xb=\fb(\xb),\quad \xb\in\R^n,
\end{equation}
where $\fb:\R^n\supset \Delta \to \R^n$ is smooth, and $x_0$ is a singularity, i.e.  $\fb(\xb_0)=0$, the celebrated Hartman-Grobman theorem \cite{Chi} states that the linearization of \eqref{sys} or equivalently the  set of the eigenvalues $\lambda_1,\ldots,\lambda_n$ of the Jacobian matrix $D\fb(\xb_0)$ characterizes the local qualitative behavior of the trajectories provided that the eigenvalues have non-zero real part, i.e. $Re(\lambda_j)\ne 0$. In this case, we say that $x_0$ is hyperbolic, otherwise we say that it is nonhyperbolic
otherwise.  To establish the local stability of nonhyperbolic singular points,  higher order terms have to be taken into account.

\smallskip

\begin{figure}
\begin{center}
\includegraphics[width=0.5\textwidth]{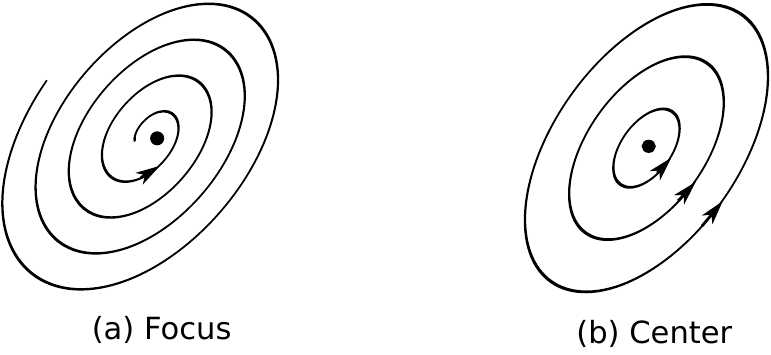}
\end{center}
\caption{An example of a  stable focus (a) and a center (b).}\label{Fig:CF}
\end{figure}

One of the simplest and well-known stability questions is the \emph{center-focus} (or \emph{center}) problem, originally defined for planar polynomial differential systems, i.e., system  \eqref{sys} when $n=2$ and $\fb$ is a system of $2$ polynomials in $\R[\xb]$ of some degree $m$. It consists of obtaining conditions on the coefficients of  $\fb(\xb)$ to distinguish between a local focus (see Fig.\,\ref{Fig:CF}(a)) or a center (see Fig.\,\ref{Fig:CF}(b)), which has been the subject of intensive research (e.g., \cite{Sib54, VulSib89, Zol94, Zol94:center, Chr94, TeiYan01, Bru06,  BriRoyYom10, RomSha, ChrLli00, Fercec2013}). 
Although the problem is open in its full generality, 
it has been solved for some important subclasses of planar 
polynomial vector fields.
As an example, consider the quadratic system defined by 
\begin{equation}\label{sys:planar:quadratic}
\begin{aligned}
\dot u&=\,\,\,\,v+a_1 u^2+a_2 uv+a_3 v^2\\
\dot v&=-u+a_4 u^2+a_5 uv+a_6 v^2,
\end{aligned}
\end{equation}
where $a_1,\ldots,a_6\in\R$.  The center conditions 
were established by Dulac \cite{Dul08} and Kapteyn~\cite{Kap2}. It is
well-known (see e.g. \cite{Zol94:center, RomSha}) that, for system
\eqref{sys:planar:quadratic}, the Bautin ideal $\B$  is generated by
the first three focus quantities of this system \cite{Bautin1952}. Moreover, the
center variety $\vv(\B)\subset\R^6$ of the ideal $\B$ has
four irreducible components, namely
\[\vv(\B)=\vv(I_{Ham})\cup\vv(I_{sym})\cup\vv(I_{\triangle})\cup\vv(I_{con}),\]
corresponding to Hamiltonian
systems, reversible systems, the Zariski  closure of those systems
having three invariant lines, and the Zariski closure of systems having an invariant
conic and an invariant cubic, respectively.  

\smallskip

The center-focus problem can also be defined for higher dimensional systems \cite{bib} and have recently been studied for a number of three-dimensional families \cite{EMRS, BGM, GMS13, Mah13:IJBC, MPS, MRS13}.
We continue this study here by applying our new symbolic-numerical
approach to a three-dimensional system
presented in Sec.~\ref{Application} with results presented
in Theorems~\ref{MainThe:01} and~\ref{MainThe:02}.

%--------------------------------------------------------------
\subsection{Computational challenges and the new approach} 
The process of solving the center-focus problem for a specific system of differential equations can be divided into three steps \cite{CL07}.  
First, the computation of certain number, say $p\in\N$, of \emph{focus quantities} (also called Lyapunov quantities), which are polynomials in the parameters of the system. Second, finding the common zeros of the polynomial system formed by the focus quantities, or more precisely the determination of the irreducible component of the variety of the ideal generated by the first $p$ focus quantities. Third, for the system restricted to each component, one checks if the necessary conditions for the existence of a center can be applied. This typically involves the application of
the Darboux theory of integrability or reduction~to~the~center~manifold.  

\smallskip

Techniques for efficient computation of Lyapunov quantities has been motivated both by mathematical and engineering problems. Over the years,
a number of algorithms have been developed \cite{Romanovski1993, Pearson1996, Gasull1997, Gasull2001, Lynch2005, Kuznetsov2008, Yu2008}. In this work, we used a method (described in~\cite{EMRS}) for computing the focus quantities for a system in dimension three, which is based on the equivalence of the existence of a center and a local analytic first integral in the neighborhood of a singular point 
(more details are provided in Sec.~\ref{sec:Lyapunov}).  
The advantage of this approach is that it allows to avoid center manifold approximation, which is especially important since its power series approximation of analytic or even polynomial systems need not~converge~(e.g.,~see~\cite{Aulbach1985, Sij85, MRS13}). 

%the center manifold entirely using the ideas developed in \cite{EMRS} (and reviewed in Section \ref{sec:Lyapunov})."

%\begin{itemize}
%\item 1885 Ponicare, Lyapunov 1892
%\item 1996 Pearson et\,al. \cite{Pearson1996}
%\item 1997 Gasull \& Prohens \cite{Gasull1997} First Lypunov
%\item 2001 Gasull \& Torre \cite{Gasull2001} Lyapunov for degenerate system
%\item 2005 Lynch \cite{Lynch2005}
%\item 2008 Kuznetsov \cite{Kuznetsov2008b} and \cite{Kuznetsov2008}
%\item 2008 Yu Chen \cite{Yu2008} Comparison of three methods. 
%\end{itemize}

\smallskip

From the computational point of view, the biggest obstacle in solving the center-focus problem for a specific system is the determination of the irreducible components of the variety (i.e., solution set) 
defined by a certain number of focus quantities. The most common approach \cite{Aziz2012, Gine2014, Fercec2014} is the application of computer algebra algorithms for computing the primary decomposition of the ideal generated by the focus quantities such as Gianni-Trager-Zacharias~(GTZ)~\cite{minGTZ} or Shimoyama-Yokoyama (SY)~\cite{SY}, which have been implemented in various symbolic packages (e.g. {\sc Singular} \cite{sing}, or {\sc Macaulay} \cite{Macaulay2}). The computational difficulty related with Gr\"obner basis calculation over the field of characteristic zero was eased by implementation of modular arithmetics \cite{Winkler1988, Edneral1997, Romanovski2011}, and successfully used in numerous problems \cite{Fercec2011, Han2012, Fercec2014, Valls2015}.  Unfortunately, in practice, the application of algorithms that use Gr\"obner bases (also with modular arithmetics) is computationally very heavy and the center conditions can only be obtained for specific systems with few parameters. In this paper, we replace this particular step and find the common zeros of the polynomial systems formed by focus quantities using numerical algebraic geometry techniques (for more details, see Sec.\,\ref{sec:NAG} and the books \cite{BertiniBook,SW05}). 
The parallelizablity of numerical algebraic geometry together
with a regeneration based approach~\cite{Regen,HW14}
and exactness recovery~\cite{RecoverExact} provides a natural
alternative to Gr\"obner basis methods.  
In particular, for the first time, we are able to solve the center-focus problem for a quadratic, three-dimensional system described next.

%Due to computational complexity, the center conditions for system \eqref{StanFor} with $h$ defined by \eqref{equaG} have been determined only after setting at least two coefficients to zero.  Although one can 
%compute the first 8 nonzero focus quantities, the standard symbolic algorithms such as Gianni-Trager-Zacharias \cite{minGTZ} or Shimoyama-Yokoyama \cite{SY}  (both implemented in various computer algebra packages such as {\sc Singular} \cite{Singular}) were not able to provide the decomposition of the Bautin ideal into primes, even for 5 parameter families. 

%\smallskip
%
%The study of the stability for multiparameter families of differential equations is computationally a very demanding task. Surprisingly, despite an intensive research \cite{CL07} even the center conditions for a planar polynomial cubic system \eqref{sys} remain unknown!  Still, recent advancements in  computer algebra allowed to solve a number of stability and bifurcation problems. Using various symbolic manipulators (e.g. Maple, Mathematica, Singular, etc.) made it possible to.....
%\marginpar{Describe Singular, Focus quantities}

%--------------------------------------------------------------
\subsection{An application}\label{Application}
Consider a third-order differential equation of the form
\begin{equation}\label{maineq}
\dddot u=\ddot u +\dot u+u+f(u,\dot u,\ddot u),
\end{equation}
where  $f=f(u,\dot u,\ddot u)\in\R[u,\dot u,\ddot u]$ is a polynomial of degree $m$. Following \cite{Mah13:IJBC}, we can equivalently write
\begin{equation}\label{StanFor}
\dot u=-v+h(u,v,w),\qquad
\dot v=u+h(u,v,w),\qquad
\dot w=-w+h(u,v,w),
%\end{aligned}
\end{equation}
where $h(u,v,w)=f(-u+w,v-w,u+w)/2$, which we call the {\it standard form} of system \eqref{maineq}. Note that the origin of  \eqref{StanFor} is a nonhyperbolic singularity at which the associated Jacobian has two purely imaginary eigenvalues $\lambda_{1,2}=\pm i$ and $\lambda_{3}=-1$. Various dynamic aspects of systems of the form \eqref{StanFor} have recently been considered, including the center conditions~\cite{BGM, DiaMel10, EMRS, MPS}, limit cycle bifurcations~\cite{WLC, MRS13}, Lie symmetries \cite{GMS13}, and isochronicity \cite{RMF13}.  In particular, the center conditions on the local center manifold for system \eqref{StanFor}, where
\begin{equation}\label{equaG}
h(u,v,w)=a_1 u^2+a_2 v^2+a_3 w^2+a_4u v+a_5 u w+a_6 v w,
\end{equation}
were studied in \cite{Mah13:IJBC}.  Although it was possible to compute the first eight focus quantities,  standard symbolic algorithms (e.g. GTZ and SY) were not able to provide the decomposition of the Bautin ideal into primes for a general six-parameter system, even over the field of non-zero characteristics. On the other hand, the application of our hybrid 
approach using numerical algebraic geometry to decompose 
described in this paper, allowed us to obtain the center conditions for a general six-parameter system \eqref{StanFor}.
\begin{theorem}\label{MainThe:01}
The system \eqref{StanFor} with $h(u,v,w)$ as in \eqref{equaG}
admits a center on the local center manifold if and only if one 
of the following holds:

\smallskip

\begin{itemize}
\item[\rm (1)] $a_1 = a_2 = a_4 = 0$
\item[\rm (2)] $a_1-a_2 = a_3 = a_5  = a_6 = 0$
\item[\rm (3)] $a_1+a_2 = a_3 = a_5 = a_6 = 0$
\item[\rm (4)] $a_1+a_2 = 2a_2-a_3+a_6 = a_3-a_4-2a_5 = 2a_4+3a_5+a_6 = 0$
\item[\rm (5)] $2a_1 - a_6 = 2a_2 + a_5 = 2a_3-a_5+a_6 = a_4+a_5+a_6 = 0$
\item[\rm (6)] $a_1-a_2 = 2a_2+a_6 = a_4 = a_5 + a_6 = 0$
\item[\rm (7)] $2a_1 + a_2 = 2a_2 + a_6 = 4a_3 + 5a_6 = a_4 = 2a_5 - a_6 = 0$.
\end{itemize}
\end{theorem}

As an easy conclusion, note that the irreducible components of the center variety (i.e. the variety of the Bautin ideal generated by the focus quantities) of system \eqref{StanFor} for quadratic $h$ \eqref{equaG} are vector subspaces of its six-dimensional parameter space, which was conjectured~in~\cite{Mah13:IJBC}.

\subsection{Outline}
The rest of the paper is organized as follows.
Section~\ref{sec:Lyapunov} summarizes focus quantities
and their computation.
Section~\ref{sec:NAG} summarizes the numerical algebraic
geometric solving approach along with exactness recovery
method used to prove Theorem~\ref{MainThe:02} in
Section~\ref{sec:Proof}.  Appendix~\ref{sec:DK}
presents the Dulac-Kapteyn criterion of quadratic planar systems
with Appendix~\ref{sec:Darboux} summarizing Darboux theory
of integrability.

%This paper aims to show how the combination of symbolic methods for computing the focus quantities (see Sec.\,\ref{sec:Lyapunov}) for a three dimensional system together with numerical algebraic geometry  solving algorithms along numeric-symbolic exactness recovery methods  (see Sec.\,\ref{sec:NAG}) can further push the boundaries related with obtaining the center conditions. Specifically, we will solve completely the center problem for system \eqref{StanFor} without any restrictions on the parameters of $h$ given by \eqref{equaG}.  

%----------------------------------------------------------
\section{Focus quantities computation in $\R^3$}\label{sec:Lyapunov}

This section is a review of the method described in \cite{EMRS} (see also \cite{Mah13:IJBC, MRS13}) for studying the center problem on a center manifold for vector fields in dimension three. Let $X: U\rightarrow \mathbb R^3$  be a real analytic vector field, such that $DX(0)$ has one non-zero and two purely imaginary eigenvalues. By an invertible linear change of coordinates and a possible rescaling of time, the system of differential equations ${\bf \dot{u} }= X({\bf u})$ can be written in the form
\begin{equation}\label{NorFor}
\begin{array}{lclcl}
\dot{u} & = & -v+P(u,v,w) \vspace{0.2cm}\\
\dot{v} & = & u+Q(u,v,w) \vspace{0.2cm}\\
\dot{w} & = & \beta w+R(u,v,w),
\end{array}
\end{equation}
where $\beta$ is a non-zero real number. Let $X=(-v+P){\p}/{\p u} + (u+Q){\p}/{\p v}+(\beta w+R){\p}/{\p w}$
denote the corresponding vector field.
A~{\it local first integral} of system \eqref{NorFor} is a nonconstant differentiable function $H$ defined in a neighborhood of the origin 
in $\mathbb R^3$ mapping into $\mathbb R$ that is constant on trajectories of \eqref{NorFor}, equivalently, $H$ satisfies
\begin{equation}
\label{IntPri}
 XH:=(-v+P)\frac{\partial H}{\partial u}+( u+Q)\frac{\partial H}{\partial v}+(\beta w+R)\frac{\partial H}{\partial w}\equiv 0
\end{equation}
sufficiently close to the origin. A {\it formal first integral} for system \eqref{NorFor} is a non-constant formal power series $H$ in $u$, $v$ and $w$ such that when $P$, $Q$, and $R$ are expanded in power series at the origin, every coefficient in the formal power series in \eqref{IntPri} is zero. 
%When $w$ and $\dot{w} $ are absent from system \eqref{NorFor} so that the system is in $\mathbb R^2$, the singularity at the origin is either a center (a punctured neighborhood is composed entirely of periodic orbits) or a focus (every trajectory near the origin spirals towards the origin, or every trajectory does so in reverse time). The center problem is the problem of distinguishing between the two cases. It was solved by Poincar\'e and Lyapunov in terms of the existence or non-existence of a local first integral. A proof appears in \cite{RomSha}.

\smallskip

It is well-known that system \eqref{NorFor} admits a local center manifold $W^c_{loc}$ at the origin, e.g.,~see~\cite[Thm.~$5.1$]{KUZ}. One of the main tools for detecting a center on a center manifold is the following theorem (see, e.g., \cite{bib, EMRS}).

\begin{theorem}  
\label{PoiLya}
The following statements are equivalent.
\begin{itemize}
\item[\rm (a)] The origin is a center for $X\mid_{W^c_{loc}}$.
\item[\rm (b)] System \eqref{NorFor} admits a local analytic first integral at the origin.
\item[\rm (c)] System \eqref{NorFor} admits a formal first integral at the origin.
\end{itemize}
\end{theorem}
In fact, a real analytic local first integral from statement (b) (as well as a formal first integral from statement (c)) can always be chosen to be of the form $H(u,v,w) = u^2 +v^2 +\cdots$ where the dots mean higher order terms 
in a neighborhood of the origin in $\mathbb R^3$.

\smallskip

The equivalence of statements (a) and (b) is called the {\it Lyapunov Center Theorem} with a proof presented in, e.g., \cite{bib}.  By this theorem, we can restrict our efforts to investigate the conditions for the existence of a first integral $H$ which is equivalent to determine necessary and sufficient conditions for the existence of a center or a focus on the local center manifold.  

\smallskip

From now on, we assume that $P$, $Q$ and $R$ in \eqref{NorFor} are polynomials. We begin by introducing the complex variable $x = u + iv$.  The first two equations in \eqref{NorFor} are equivalent to a single equation $\dot{x} = ix + \cdots$, where the dots represent a sum of homogeneous polynomials of degrees between $2$ and $n$. 
Let $\bar{x}$ denote the complex conjugate of $x$.
We add to this equation its complex conjugate, replacing $\bar{x}$
everywhere by $y$ which is regarded as an independent complex variable
and replacing $w$ by $z$ simply as a notational convenience.
This yields the following complexification of \eqref{NorFor}:
\begin{equation}
\label{NorFor2}
\begin{array}{lcl}
\dot{x} & = & ix+\displaystyle\sum^{n}_{p+q+r=2}a_{pqr} x^p y^q z^r, \vspace{0.2cm}\\
\dot{y} & = & -iy+\displaystyle\sum^{n}_{p+q+r=2}b_{pqr} x^p y^q z^r,  \vspace{0.2cm}\\
\dot{z} & = &\beta z+\displaystyle\sum^{n}_{p+q+r=2}c_{pqr} x^p y^q z^r, 
\end{array}
\end{equation}
where $b_{qpr} = \bar{a}_{pqr}$ and $c_{pqr}$ are such that $\sum^{n}_{p+q+r=2} c_{pqr} x^p \bar{x}^q w^r$ is real for all $x\in \mathbb C$ and $w\in\mathbb R$. Let $X$ be the corresponding vector field of system \eqref{NorFor2} on $\mathbb C^3$.  Existence of a first integral $H(u,v,w) = u^2 + v^2 + \cdots$ for system \eqref{NorFor} is equivalent to the existence of a first integral for system \eqref{NorFor2}, denoted again by $H$, of the form
\begin{equation}
\label{FirstInt}
H(x,y,z)=xy+\sum_{j+k+\ell=3} v_{jkl} x^j y^k z^\ell.
\end{equation}

\smallskip

We now investigate the existence of a first integral $H$ for system \eqref{NorFor2} by computing the coefficients of $XH$ and equating them to zero. When $H$ has the form  \eqref{FirstInt}, 
the coefficient $g_{jk\ell}$ of $x^{j} y^{k} z^{\ell}$ in $XH$ can be calculated explicitly (see \cite{EMRS}). 
Except when $j = k$ and $\ell = 0$, the equation $g_{jk\ell}$ = 0 
can be solved uniquely for $\nu_{jk\ell}$ in terms of the known 
quantities $\nu_{\alpha\beta\gamma}$ with $\alpha+\beta+\gamma<j+k+\ell$. 
A formal first integral $H$ thus exists if $g_{KK0} = 0$ for all $K\in\mathbb  N$. An obstruction to the existence of the formal series $H$ occurs when the coefficient $g_{KK0}$ is non-zero.  This coefficient is the $K^{\rm th}$ {\it focus quantity} and it can be expressed as
\begin{equation}\label{e:focus.quantity}
{\small
\begin{aligned}
g_{KK0} &=\sum_{\substack{j+k = 2 \\ j \ge 0, k \ge 0}}^{ 2K - 1 }
           \left(
           j \, a_{K-j+1, K-k,0} + k \, b_{K-j, K-k+1,0}
           \right)
           v_{j,k,0}~~+~\sum_{\substack{j+k = 2 \\ j \ge 0, k \ge 0}}^{ 2K - 2 }
           c_{K-j, K-k,0} \, v_{j,k,1},
\end{aligned}
}
\end{equation}
where we have made the natural assignments $v_{110} = 1$ and $v_{\alpha \beta \gamma} = 0$ for 
$\alpha + \beta + \gamma = 2$ but $(\alpha, \beta, \gamma) \ne (1,1,0)$. 
We know $g_{110} = 0$ and $g_{220}$ is uniquely determined, but the remaining ones depend on the choices made for $v_{KK0}$, $K \in \mathbb N$, $K\geq 2$. Once such an assignment is made, $H$ is determined and satisfies
\[
XH(x,y,z) = g_{220} (xy)^2 + g_{330} (xy)^3 +\cdots .
\]
It is known that if at least one focus quantity is non-zero
for a choice of $v_{KK0}$, then the same is true for every other choice of the $v_{KK0}$. The vanishing of all focus quantities, i.e., 
\begin{equation}\label{fq:vanish}
g_{KK0} = 0\qquad \text{for}\qquad K\geq2
\end{equation}
is both a necessary and sufficient condition for the existence of a center on the center manifold, otherwise there is a focus (see \cite{EMRS}). 

\smallskip

By Hilbert's basis theorem, there exists $K_0\geq 2$ such that the set of solutions of $g_{KK0} = 0$ for all $2\leq K\leq K_0$ is equivalent set defined by an infinite system  \eqref{fq:vanish}.
Since such a $K_0$ is not known {\em a priori},  we will apply an iterative approach that solves $g_{KK0} = 0$ for $2\leq K\leq M+1$ given the solution set of $g_{KK0} = 0$ for $2\leq K\leq M$.  Without knowing $K_0$,  solving using any $M\geq2$ does always yield necessary conditions.

%================================
\section{Numerical algebraic geometry}\label{sec:NAG}

Symbolic methods, such as Gr\"obner basis techniques, take an algebraic
viewpoint for solving systems of polynomial equations.  
In broad terms, they manipulate equations to obtain new relations
describing the solution set.  An alternative approach is 
to use a geometric viewpoint which manipulates solution sets.
Following a numerical algebraic geometry approach, solution sets
are represented by {\em witness sets} that we discuss below.
A more detailed comparison of symbolic and numerical approaches 
is provided in \cite{Bates2014}.

\smallskip

The field of numerical algebraic geometry grew out of 
the use of homotopy continuation for computing isolated solutions
to a system of polynomial equations.  We will first briefly
explain using basic homotopy continuation on a polynomial 
system $F:\bC^N\rightarrow\bC^N$, that is, $F(x) = 0$ defines a system
of $N$ polynomial equations in $N$ variables.  
The idea is to select another polynomial system $G:\bC^N\rightarrow\bC^N$
related to $F$ such that $G(x) = 0$ is ``easy'' is to solve.  
The simplest example is $G_i(x) = x_i^{d_i} - 1$ where $d_i = \deg F_i$,
but there is a wide range of constructing the so-called {\em start systems} 
which exploit structure in $F$ (see \cite{BertiniBook,SW05} for a broad overview).
Let $S\subset\bC^N$ denote the set of isolated nonsingular solutions of $G = 0$.

\smallskip

The next step is to construct a homotopy $H:\bC^N\times\bC\rightarrow\bC^N$
connecting $G$ and $F$, say 
$$H(x,t) = F(x)\cdot(1-t) + \gamma\cdot t\cdot G(x)$$
for a randomly selected $\gamma\in\bC$.  For each $s\in S$, the homotopy $H$ 
defines a solution path~$x_s(t)$ such that $x_s(1) = s$ and $H(x_s(t),t)\equiv 0$
with the goal of computing the endpoint $x_s(0) = \lim_{t\rightarrow0^+} x_s(t)$.  
In fact, this limit is either a point in $\bC^N$ 
which must be a solution of $F = 0$ or the path is said to be diverging to infinity.
By differentiating $H(x_s(t),t)$ with respect to $t$, one obtains the 
 Davidenko differential equation
$$J_x H(x_s(t),t)\cdot \dot{x}_s(t) = -J_t H(x_s(t),t).$$
By including the randomly selected $\gamma$, called the ``gamma trick,''
the Jacobian matrix $J_x H$ is invertible along the path for $t\in(0,1]$
with probability one and thus one can use predictor-correct techniques
to track the solution path $x_s(t)$ starting at $x_s(1) = s$ 
in order to approximate $x_s(0)$.  We refer the interested reader
to \cite{BertiniBook,SW05} for more details about path tracking
and using endgames to estimate $x_s(0)$.
In the end, the set $E\subset\bC^N$ of convergent endpoints of
all the paths $x_s(t)$ for $s\in S$ is a superset 
of the isolated nonsingular solutions of $F = 0$.  

\smallskip

We now turn our attention to computing the solution set of $F = 0$, denoted $\Var(F)\subset\bC^N$, for a polynomial system $F:\bC^N\rightarrow\bC^n$.  
Geometrically, $\Var(F)$ can be decomposed into a union of irreducible components
$\Var(F) = \cup_{i=1}^r V_i$.  This corresponds algebraically to a prime decomposition
of the radical ideal generated by $F$, namely $\sqrt{I(F)} = \cap_{i=1}^r I(V_i)$.
Numerical algebraic geometry describes an irreducible decomposition of $\Var(F)$
by computing a witness set for each $V_i$, called a numerical irreducible decomposition.

\smallskip

Suppose that $V$ is an irreducible component of $\Var(F)$ for some polynomial system $F$.
A witness set for $V$ is the triple $\{F,\sL,W\}$ where $\sL\subset\bC^N$
is general linear subspace of codimension $d = \dim V$ and 
$W = V\cap\Var(\sL)$ so that $|W| = \deg V$.  Here, the definition of 
general means that $\sL$ intersects $V$ transversely, which is a Zariski
open condition on the Grassmannian of codimension $d$ linear subspaces in $\bC^N$.
The books \cite{BertiniBook,SW05} provide for more information about 
witness sets including performing computations on irreducible components
which have multiplicity $> 1$ with respect to $F$.

\smallskip

A witness set for an irreducible component $V\subset\bC^N$ 
facilitates additional computations that can be performed on $V$.  
Of particular interest to the problems discussed in this article
include 
the recovery of exact polynomials that vanish on $V$,
determining the existence of real points in $V$, 
and intersecting $V$ with another solution set.

\smallskip

With an input polynomial system with exact coefficients, e.g., in $\Q$,
one often would like exact output.  Although the internal computations
and witness sets rely upon numerical approximations, there exist techniques
for recovering exact answers which can then be 
verified using exact symbolic methods, which is typically computationally inexpensive.
For the problems at hand here, we
use the exactness recovery technique described in \cite{RecoverExact}
which uses a sufficiently accurate numerical approximation of
a sufficiently general point on $V$ to compute polynomials
with integer coefficients that vanish on $V$.  
This method is based on using a lattice-base reduction technique
such as LLL \cite{LLL} or PSLQ \cite{PSLQ}.

\smallskip

In many applications, only real solutions or components which contain
real points are of interest, which is the case here.  
The approach of \cite{RealPoints} uses critical points conditions
of the distance function to determine if $V$, represented by a witness set,
contains real points.  If~$V\cap\R^N = \emptyset$, then we can disregard 
this component from further computations.

\smallskip

We conclude this section by describing the intersection approach
built from witness sets which is used in the subsequent section.
For this situation, we consider a sequence of polynomial systems
of interest, namely $F_k = \{f_1,\dots,f_k\}$ for $k\geq 1$.
Given witness sets for the irreducible components of $\Var(F_k)$,
our goal is to compute witness sets for the irreducible components
of $\Var(F_{k+1}) = \Var(F_k)\cap\Var(f_{k+1})$.
For consistency, we assume $\Var(F_k)\subset\bC^N$ 
since one can easily adjust the methods to work on projective space
which will arise below since each $g_{KK0}$ is homogeneous in $a_1,\dots,a_6$,
i.e., $\Var(g_{KK0})$ is naturally a hypersurface in $\bP^5$.

\smallskip

For the base case, we need to decompose the hypersurface $\Var(F_1) = \Var(f_1)$,
which can be readily performed,~e.g.,~via~\cite{Hypersurfaces}.

\smallskip

Now, suppose that we are given witness sets for 
the irreducible components $V_{k,1},\dots,V_{k,n_k}$ 
of $\Var(F_k)$.  For each $j\in\{1,\dots,n_k\}$, we need to compute
$V_{k,j}\cap\Var(f_{k+1})$ using the provided witness set for~$V_{k,j}$,
say $\{F_k,\sL_{k,j},W_{k,j}\}$.
Clearly, if $f_{k+1}$ vanishes identically on $V_{k,j}$,
we know that $V_{k,j}$ is an irreducible component of $\Var(F_{k+1})$.
Thus, we shall assume that $V_{k,j}$ is not contained in the hypersurface
$\Var(f_{k+1})$ so that $V_{k,j}\cap\Var(f_{k+1})$ 
is either empty or consists of irreducible components of 
dimension one less than $V_{k,j}$.
With this assumption, if $d := \dim V_{k,j}$ is zero, we know
$V_{k,j}\cap\Var(f_{k+1}) = \emptyset$.  Thus,
we also assume that $d > 0$.  

\smallskip

We compute $V_{k,j}\cap\Var(f_{k+1})$ using a 
regenerative intersection approach developed in \cite{HW13,HW14}
which builds on the diagonal intersection \cite{Diagonal}
and the regenerative cascade \cite{RegenCascade,Regen}.
It can be performed using {\sc Bertini}~\cite{Bertini}.
To perform this computation, we select a general hyperplane~$\sH$
and codimension $d-1$ linear space $\sK$ so that $\sL_{k,j} = \sH\cap\sK$.
Our first goal is to compute the finite set of points 
$V_{k,j}\cap\sK\cap\Var(f_{k+1})$ given $W_{k,j} = V_{k,j}\cap\sK\cap\sH$.
If $g = \deg f_{k+1}$, we select general hyperplanes $\sH_1,\dots,\sH_g$
and compute $V_{k,j}\cap\sK\cap\sH_\ell$ for $\ell = 1,\dots,g$ by standard
homotopy continuation from $V_{k,j}\cap\sK\cap\sH$.  
Thus, we have computed
$$V_{k,j}\cap\sK\cap(\cup_{\ell=1}^g \sH_\ell)$$
which, again by standard homotopy continuation, can be used to compute
$V_{k,j}\cap\sK\cap\Var(f_{k+1})$.  

\smallskip

Now, to compute witness sets for the irreducible components 
of $V_{k,j}\cap\Var(f_{k+1})$, 
which will have the form $\{F_{k+1},\sK,\bullet\}$,
we simply need to partition
the set of points $V_{k,j}\cap\sK\cap\Var(f_{k+1})$ 
into subsets corresponding to distinct irreducible components.  
This is accomplished by using random monodromy loops \cite{Monodromy} 
with the decomposition certified by the trace test \cite{Trace}.
In short, the idea is to move the linear space $\sK$ in a random loop
in the Grassmannian and observe which points are connected by such paths:
the points in $V_{k,j}\cap\sK\cap\Var(f_{k+1})$ 
which are path connected by such random loops 
must lie on the same irreducible component. 
Thus, monodromy loops yield necessary conditions.  
The (linear) trace test yields a sufficient condition which follows
from the fact that, as the linear space $\sK$ is moved in parallel,
the centroid of points arising from a union of irreducible components
must move linearly.  
One can read about further details of computing 
such decompositions in \cite{BertiniBook,SW05}.

\section{Center conditions for a three dimensional quadratic system}\label{sec:Proof}
The following provides a proof of Theorem \ref{MainThe:01}. 
We note that, without loss of generality, we can always assume that either $a_6=0$ or $a_6=1$.
The latter follows immediately by the change of variables $(u,v,w)\mapsto (x/a_6, y/a_6, z/a_6)$ and rescalling of time  $dt=a_6 d\tau$. 
Thus, the seven cases in Theorem~\ref{MainThe:01} 
can be split into ten cases, five each for $a_6 = 0$
and $a_6 = 1$.  After showing these ten cases, we then describe
how the seven cases of Theorem~\ref{MainThe:01} follow.

\begin{theorem}\label{MainThe:02}
Consider system \eqref{StanFor} with $h(u,v,w)$ as in \eqref{equaG}.  

\smallskip

The system \eqref{StanFor} with $a_6 = 0$ admits a center on the local center manifold if and only if one of the following holds:
\begin{itemize}
\item[\rm (a)]    $a_1-a_2=a_3=a_5=0$						% l1-s,\quad ($\R$-plane)    
\item[\rm (b)]   $a_1+a_2=a_3=a_5=0$  					% l2-s,\quad ($\R$-plane)  
\item[\rm (c)]   $a_1=a_2=a_4=0$;			% p1-s,\quad ($\R$-plane)
\item[\rm (d)]    $a_1+a_2=2a_1+a_3=6a_1-a_4=4a_1 + a_5=0$	% l3-s,\quad ($\R$-line)    
\item[\rm (e)]    $a_1=a_2+a_3=2a_2-a_4=2a_2+a_5=0$. % l4-s,\quad ($\R$-line)
\end{itemize}

The system \eqref{StanFor} with $a_6 = 1$ admits a center on the local center manifold if and only if one of the following holds:
\begin{itemize}
\item[\rm (f)]    $a_1=a_2=a_4=0$ 				% p1,\quad ($\R$-plane)
\item[\rm (g)]    $2a_1-1=a_4+a_5+1=2a_2+a_5=2a_3-a_5+1=0$ % l1,\quad ($\R$-line)
\item[\rm (h)]    $2a_1+1=2a_2+1=a_4=a_5+1=0$	% l2,\quad ($\R$-line)
\item[\rm (i)]    $a_1+a_2=4a_2-a_5+3=6a_2+a_4+5=2a_2-a_3+1=0$ % l3,\quad ($\R$-line)
\item[\rm (j)]    $4a_1-1=2a_2+1=4a_3+5=a_4=2a_5-1=0$.	% s1,\quad ($\R$-point)  
\end{itemize}
\end{theorem}

\medskip

\noindent {\bf Necessary conditions}.

We first consider $a_6 = 0$ and take $(a_1,\dots,a_5)\in\bP^4$.  
Using the notation from Section~\ref{sec:NAG},
$\Var(F_2)$ and $\Var(F_3)$ are irreducible of codimension $1$ and $2$
of degree $2$ and $8$, respectively.  Now, $\Var(F_4)$
has codimension $3$ and decomposes into the following components:
$5$ linear spaces, $3$ of
multiplicity $1$ and $2$ of multiplicity $3$, and
an irreducible algebraic set of degree~$39$.  
The three linear spaces of multiplicity $1$ are (a), (b), and (c).  
The other two linear spaces are complex conjugates of each other with
their union is defined in $\bP^4$ by
$$a_1 + a_2 = 4a_2^2 + a_4^2 = a_5 = 0.$$
Since the real points on this union are contained in (c),
we only need to further investigate the degree $39$ component, denoted 
$X_{4,6}$, which is not contained in $\Var(g_{550})$.  Regenerating from $X_{4,6}$ 
to compute $X_{4,6}\cap\Var(g_{550})$ yields 
$189$ distinct points in $\bP^4$, of which $19$ correspond to real points.
There are $14$ real points that do not lie on (a), (b), or (c)
of which only $2$ satisfy $g_{660}=0$, namely (d) and (e).  
We note that (e) has multiplicity $2$ with respect to $F_5$.

We next consider $a_6 = 1$ and take $(a_1,\dots,a_5)\in\bC^5$.  
Similar to the case above, 
$\Var(F_2)$ and $\Var(F_3)$ are irreducible of codimension $1$ and $2$
of degree $2$ and $8$, respectively.
Also, $\Var(F_4)$
has codimension $3$ and decomposes into the following components: 
$3$ linear spaces, one having 
multiplicity $1$, namely (f), with the other $2$ having multiplicity $3$, and
an irreducible algebraic set of degree $41$.  
As above, the two linear spaces of multiplicity $3$ are complex conjugates
of each other with their union defined in $\bC^5$ by 
$$a_1 + a_2 = 4a_2^2 + a_4^2 = 2a_2 + a_4 a_5 = 2a_2 a_5 - a_4 = a_5^2 + 1 = 0.$$
Since there are no real points on this union, we only need to further investigate
the degree~$41$ components, denoted $X_{4,4}$, which is not contained in
$\Var(g_{550})$.  Regenerating $X_{4,4}$ yields $4$ irreducible components 
of $X_{4,4}\cap\Var(g_{550})$ not contained in (f) or the hyperplane $a_5^2 + 1 = 0$.
Three of these are the lines (g), (h), and (i) with the fourth
being an irreducible curve of degree $244$, denoted $X_{5,4}$, 
not contained in $\Var(g_{660})$.  Regenerating $X_{5,4}$
yields $71$ distinct real points not contained in the 
hyperplane $a_5^2 + 1 = 0$ nor satisfying (f), (g), (h), or (i).
Of these, only one satisfies $g_{770} = 0$, namely (j).

\noindent {\bf Sufficient conditions}. 

\smallskip

{\it Cases (a) and (b).} If the condition (a)  (resp. (b)) holds, system \eqref{StanFor} reduces to
\[
\begin{array}{lcr}
\dot{u} & = &-v+ a_1 u^2+a_2 v^2+a_4 uv,\\
\dot{v} & = &u + a_1 u^2+a_2 v^2+a_4 uv,\\
\dot{w} & = &-w+ a_1 u^2+a_2 v^2+a_4 uv,
\end{array}
\]
with $a_2= a_1$ (resp. $a_2=- a_1$). Note that by Theorem \ref{PoiLya}, it is enough to show that this system admits a local analytic first integral at the origin. Since the first two equations are decoupled from the third we only need to show that
\begin{equation}
\label{sys:01}
\dot{u}  = -v+ a_1 u^2+a_2 v^2+a_4 uv, \quad \dot{v}  = u + a_1 u^2+a_2 v^2+a_4 uv,
\end{equation}
admits a local analytic first integral.  In fact, if $a_4\neq 0$ and $a_2=a_1$, system \eqref{sys:01} has the inverse integrating factor 
\[
V(u,v)=-a_4+a_4 \left( a_4+2 a_1 \right)  \left( x-y
 \right) +a_1\left( a_4+2 a_1 \right) ^{2} \left( {
x}^{2}+yx+{y}^{2} \right).
\]
As $V(0,0)=-a_4$, it follows that system \eqref{sys:01} has a first integral defined at the origin. If $a_4=0$ and $a_2=a_1$, applying Theorem \ref{QuaCen}(ii) with $a=c=a_1$, $b=d=-a_1$, $A=2a_1$ and $B=-2a_1$, we have that \eqref{sys:01} has a center at the origin and so it is integrable.  The case $a_2=-a_1$ (i.e. case (b)) is analogous, since
\[
V(u,v)=1+ \left( 2 a_1-a_4 \right) x+ \left( 2 a_1+ a_4
 \right) y-a_1 a_4 x^{2}-{a_4}^{2}xy+a_1a_4 {y}^{2},
\]
is an inverse integrating factor for system \eqref{sys:01}, which is also nonzero at the origin. 

\smallskip

{\it Case (c)}. In this case system \eqref{StanFor} becomes
\[
\begin{array}{lcr}
\dot{u} & = &-v+ a_3 w^{2}+a_5 uw,\\
\dot{v} & = &u+ a_3 w^{2}+a_5 uw,\\
\dot{w} & = &-w+ a_3 w^{2}+a_5 uw.
\end{array}
\]
Note that w = 0 is invariant and is a center manifold for this system. Moreover, the restriction of the associated vector field to w = 0 gives rise to a linear center.

\smallskip

{\it Case (d).} For  $a_2=-a_1$, $a_3=-2a_1$, $a_4=6a_1$ and $a_5=-4a_1$ the vector field associate to system  \eqref{StanFor} has the invariant algebraic surface $F(u,v,w)=w+a_1(u-v)^2-2a_1(v-w)^2=0$ with cofactor $K(u,v,w)=-1$. Since $F = 0$ is tangent to $w = 0$ at the origin, it is a center manifold for this system. To determine the dynamics on it first we use the change of coordinates $(u, v, w)\mapsto (x+z,y+z,z)$ that transforms  the system into 
\begin{equation}
\label{sys:03}
\begin{array}{lcl}
\dot{x} & = &-y,\\
\dot{y} & = & x+2z,\\
\dot{z} & = &-z+a_1 {x}^{2}+6 a_1 xy+4 a_1 xz-a_1 {y}^{2}+4 a_1 yz.
\end{array}
\end{equation}
The center manifold $F = 0$ in the new variables is given by $F(x,y,z)=z+a_1(x-y)^2-2a_1y^2=0$.  The restriction of system \eqref{sys:03} to $F = 0$ is given by
\[
\dot{x}  = -y, \quad \dot{v}  = x -2a_1 x^2+4a_1xy+2a_1y^2.
\]
Since this system has the following inverse integrating factor (nonzero at the origin)
\[
V(u,v)=1-4a_1\left(x-y\right) + 4a_1^2\left(x^{2}-{2}xy- {y}^{2}\right)
\]
thus  in this case system \eqref{StanFor} has a center on the center manifold. 

\smallskip

{\it Case (e).} For $a_1=0$, $a_3=-a_2$, $a_4=2a_2$ and $a_5=-2a_2$ system  \eqref{StanFor} has the invariant algebraic surface $F(u,v,w)=w-a_2(y-z)^2=0$ with cofactor $K(u,v,w)=-1$. Since $F = 0$ is tangent to $w = 0$ at the origin, it is a center manifold for this system. To determine the dynamics on it first we use the change of coordinates $(u, v, w)\mapsto (x,y+z,z)$  that transforms  the system into
\begin{equation}
\label{sys:02}
\begin{array}{lcl}
\dot{x} & = &-y-z+a_2 y^{2}+2 a_2 xy+2 a_2 yz,\\
\dot{y} & = & x+z,\\
\dot{z} & = &-z+a_2 y^2+2 a_2 yz+2 a_2xy.
\end{array}
\end{equation}
The center manifold $F = 0$ in the new variables writes as $F(x,y,z)=z-a_2y^2=0$.  The restriction of system \eqref{sys:02} to $F = 0$ is 
\[
\dot{x}  = -y+2a_2 xy+ 2a_2^2 y^3, \quad \dot{y}  = x +a_2 y^2.
\]
This system is invariant by the change of variables $(x,y,t)\mapsto (x,-y,-t)$ so that it has a center at the origin.  Hence, system \eqref{StanFor} restricted to (e) has a center on the~center~manifold. 

\smallskip

{\it Case (f).} In this case system \eqref{StanFor} becomes
\[
\begin{array}{lcr}
\dot{u} & = &-v+ a_3 w^{2}+a_5 uw+ vw,\\
\dot{v} & = &u+ a_3 w^{2}+a_5 uw+ vw,\\
\dot{w} & = &-w+ a_3 w^{2}+a_5 uw+ vw.
\end{array}
\]
It is clear that the plane w = 0 is invariant and is a center manifold for this system. Moreover, the restriction of the associated vector field to w = 0 gives rise to a linear center.

\smallskip

{\it Case (g).} If  $a_1=1/2$, $a_3=-a_2-1/2$, $a_4=2a_2-1$ and $a_5=-2a_2$, then system  \eqref{StanFor} has the invariant algebraic surface $F(u,v,w)=-2w+(u-w)^2+2a_2(v-w)^2=0$ with cofactor $K(u,v,w)=-1$. Since $F = 0$ is tangent to $w = 0$ at the origin, it is a center manifold for this system. To determine the dynamics on it first we use the change of coordinates $(u, v, w)\mapsto (x+z,y+z,z)$ that transforms  system \eqref{StanFor} with conditions (l1) into
\begin{equation}
\label{sys:07}
\begin{array}{lcl}
\dot{x} & = &-y,\\
\dot{y} & = & x+2z,\\
\dot{z} & = &-z+{x}^{2}/2+(2 a_2-1) xy+ a_2y^2+4a_2yz.
\end{array}
\end{equation}
The center manifold $F = 0$ in the new variables writes as $F(x,y,z)=-2z+x^2+2a_2y^2=0$.  The restriction of system \eqref{sys:07} to $F = 0$ is 
\[
\dot{x}  = -y, \quad \dot{v}  = x +x^2+2a_2y^2.
\]
As this system is invariant under  $(x,y,t)\mapsto (x,-y,-t)$, it follows that it has a center at the origin, i.e.  system \eqref{StanFor} under the  conditions (g) has a center on the center manifold.

\smallskip

{\it Case (h).} If  $a_1=-1/2$, $a_2=-1/2$, $a_4=0$ and $a_5=-1$. Then the vector field associate to system  \eqref{StanFor} has the invariant algebraic surface $F(u,v,w)=w+\big[ (u+w)^2+\left( v-w \right)^2  \big]/2-w^2\left( 1+a_3\right)=0$ with the cofactor $K(u,$ $v,w)=-1-2u+2a_3w$. Since $F = 0$ is tangent to $w = 0$ at the origin, it is a center manifold for this system. To determine the dynamics on it first we use the change of coordinates $(u, v, w)\mapsto (x-z,y+z,z)$, that transforms  system \eqref{StanFor} with condition (l2) into
\begin{equation}
\label{sys:06}
\begin{array}{lcl}
\dot{x} & = &-y-2z+2(1+a_3)z^2-x^2-y^2,\\
\dot{y} & = & x,\\
\dot{z} & = &-z+(1+a_3)z^2-x^2/2-y^2/2.
\end{array}
\end{equation}
The center manifold $F = 0$ in the new variables is given by 
\[
F(x,y,z)=z-(1+a_3)z^2+\dfrac{1}{2}x^2+\dfrac{1}{2}y^2=0.
\]  
The restriction of system \eqref{sys:06} to $F = 0$  gives rise to a linear center.

\smallskip

{\it Case (i).}  For  $a_2=-a_1$, $a_3=-2a_1+1$, $a_4=6a_1-5$ and $a_5=-4a_1+3$ system  \eqref{StanFor} admits an invariant algebraic surface $F(u,v,w)=w+(a_1-1)(u-w)^2+(1-2a_1)(u-w)(v-w)+(1-a_1)(v-w)^2=0$ with cofactor $K(u,v,w)=-1$. Since $F = 0$ is tangent to $w = 0$ at the origin, it is a center manifold for this system. The change of coordinates $(u, v, w)\mapsto (x+z,y+z,z)$ transforms  system \eqref{StanFor} under the conditions (i) into
\begin{equation}\label{sys:08}
\begin{array}{lcl}
\dot{x} & = &-y,\\
\dot{y} & = & x+2z,\\
\dot{z} & = &-z+a_1 {x}^{2}+(6 a_1-5) xy+2(2 a_1-1)xz-a_1y^2+4(a_1-1)yz.
\end{array}
\end{equation}
Again, in the new variables  the center manifold is given by $F(x,y,z)=z+(a_1-1)x^2+(1-2a_1)xy+(1-a_1)y^2=0$ and the restriction of \eqref{sys:08} to $F = 0$ reduces to
\[
\dot{x}  = -y, \quad \dot{v}  = x +2(1-a_1)x^2+2(2a_1-1)xy+2(a_1-1)y^2.
\]
This system has the following inverse integrating factor (nonzero at the origin)
\[
\begin{array}{lcl}
V(u,v) & = & 1+4(1-a_1)x+2\left(2a_1-1\right)y + 4(a_1-1)^2x^{2} \\ & & -4(a_1-1)(2a_1-1)xy-4(a_1-1)^2 {y}^{2}.
\end{array}
\]
Hence system \eqref{StanFor} has a center on the center manifold.

{\it Case (j).} For $a_1=1/4$, $a_2=-1/2$, $a_3=-5/4$, $a_4=0$ and $a_5=1/2$ the vector field associated to system  \eqref{StanFor} admits a polynomial first integral 
\[
\begin{array}{lcl}
H(x,y,z) & = & {x}^{2}+{y}^{2}-\dfrac{1}{2}{x}^{3}-\dfrac{1}{2}{x}^{2}y+2{x}^{2}z-\dfrac{3}{2}x{z}^{2}-
{y}^{2}x \\ \\
&  &+{y}^{2}z-\dfrac{1}{2}y{z}^{2}-\dfrac{1}{2}{x}^{
3}z+ \dfrac{5}{4}{x}^{2}{z}^{2}+\dfrac{1}{2}{y}^{2}{x}^{2}\\ \\
& &-\dfrac{3}{2}x{z}^{3}+\dfrac{1}{2}{y}^{2}{z}^{2
}  -y{z}^{3}+xyz+\dfrac{1}{8}{x}^{4}
 \\ \\
& & -{x}^{2}yz-{y}^{2}xz+2yx{z}^{2}+\dfrac{5}{8}{z}^{4},
\end{array}
\]
and so it  has a center on the center manifold.

\hfill $\Box$

\noindent {\bf Proof of Theorem~\ref{MainThe:01}.}

{\it Case (1).}  Follows from Cases (c) and (f) of Theorem~\ref{MainThe:02}.

{\it Case (2).}  Follows from Case (a) of Theorem~\ref{MainThe:02}.

{\it Case (3).}  Follows from Case (b) of Theorem~\ref{MainThe:02}.

{\it Case (4).}  Follows from Cases (d) and (i) of Theorem~\ref{MainThe:02}.

{\it Case (5).}  Follows from Cases (e) and (g) of Theorem~\ref{MainThe:02}.

{\it Case (6).}  Follows from Cases (c) and (h) of Theorem~\ref{MainThe:02}.

{\it Case (7).}  Follows from Cases (c) and (j) of Theorem~\ref{MainThe:02}.

\hfill $\Box$

\appendix

%-----------------------------------------------------------
\section{Dulac-Kapteyn criterion}\label{sec:DK}

The following theorem provides  a criterion in order to determine when a quadratic planar polynomial system has a center at the origin. It was first proven by Dulac \cite{Dul08} and Kapteyn~\cite{Kap2}, but we present the version given in \cite{Cop66}.
\begin{theorem}[Quadratic Center]
\label{QuaCen}
The system
\[
\begin{array}{lcl}
\dot{u} & = & -v-b u^2-(B+2c) uv-d v^2,\\
\dot{v} & = & u+a u^2+(A+2b) uv+c v^2,
\end{array}
\]
has a center at the origin if and only if at least one of the following three hold:
\begin{itemize}
\item[\rm (i)]  $a+c=b+d$;
\item[\rm (ii)] $A(a+c) = B(b+d)$ and $aA^3 -(3b+A)A^2B+(3c + B)AB^2 - dB^3 = 0$;
\item[\rm(iii)] $A+5b+5d=B+5a+5c=ac+bd+2a^2+2d^2 = 0$.
\end{itemize}
\end{theorem}

%-----------------------------------------------------------
\section{Basic Darboux theory of integrability}\label{sec:Darboux}
Since, by Poincar\'e theorem, the integrability is closely related to the existence of a center on a center manifold (also on the plane), we provide a short overview of the basic notions of the Darboux the	ory of integrability used in Section~\ref{sec:Proof}; for more information see \cite{Llibre2000, Goriely2001} and some applications see \cite{Gasull2002, Llibre2013, Mahdi2014}.

\smallskip

We say that  $F=F (x, y, z)\in\C[x,y,z]$ is a \emph{Darboux polynomial} and  $F=0 $ is an {\it invariant algebraic surface} of the  vector field $X$ if and only if there exists a polynomial $K (x, y, z)\in\C[x,y,z]$, the {\it cofactor} of $F$, such that $XF = K F$. A the heart of the Darboux theory of integrability is the following result \cite{Dar1878}: if there exists some number $n$ of pairs $(F_j,K_j)$ for which there exists a nontrivial dependency relation $\sum \alpha_j K_j =0$ then $F_1^{\alpha_1}\cdots F_n^{\alpha_n}$ is a first integral of $X$.

\smallskip

Consider now the planar system
\begin{equation}\label{eq:101} 
\dot{x}=P(x,y), \quad \dot{y}=Q(x,y),
\end{equation}
where $P, Q\in\R[x,y]$, and the associate vector field $X=P\p/\p x+Q\p/\p y$. Let $U$ be an open subset of $\mathbb R^{2}$, and let $R, V:U\rightarrow \mathbb R$  be two analytic functions which are not identically zero on $U$. We say that $R$ is an {\it integrating factor} of this polynomial system on $U$ if one of the following three
equivalent conditions holds
\[
\displaystyle \frac{\partial
RP}{\partial x}=-\frac{\partial RQ}{\partial x},\quad {\rm
div}(RP,RQ)=0,\quad XR=-R\; {\rm div}(P,Q),
\]
where ${\rm div}$ denotes the divergence. The first integral $H$ associated to the integrating factor $R$ can be easily obtained by
\[ 
H(x,y)=\int R(x,y)P(x,y)dy+h(x),
\]
where $h(x)$ is chosen such that it satisfies $\partial H /\partial x=-RQ$. Note that
$\partial H/\partial y=RP$, so that $XH\equiv 0$. The function $V$ is an {\it inverse integrating}
factor of the polynomial system (\ref{eq:101}) on $U$ if
\begin{equation}
\label{eq:102} P\frac{\partial V}{\partial x}+Q\frac{\partial
V}{\partial y}=\left(\frac{\partial P}{\partial x}+\frac{\partial
Q}{\partial y}\right)V.
\end{equation}
We note that $\{V=0\}$ is formed by orbits of system
(\ref{eq:101}) and $R= 1/V$ defines on $U\setminus \{V=0\}$ an
integrating factor of (\ref{eq:101}). We note that if $P$ and $Q$ are quadratic polynomials and the origin of system \eqref{eq:101}  is a center, then there always exits a polynomial function $V:\mathbb R^2\rightarrow \mathbb R$ of degree $3$ or $5$ satisfying equation (\ref{eq:102}), see \cite{LliMahFer07}.

\section*{Acknowledgments}

\noindent C.P. was partially supported by Program CAPES/DGU Process 8333/13-0 and by FAPESP-Brazil Project 2011/13152-8.  J.D.H. was partially 
supported by NSF DMS-1262428 and ACI-1460032, Sloan Fellowship, and DARPA YFA.

\bibliography{bibliografia}
\bibliographystyle{amsplain}

\end{document}